\documentclass{scrartcl}

\usepackage[utf8]{inputenc} 
\usepackage[T1]{fontenc}  
\usepackage[english]{babel}
\usepackage{amsmath} 
\usepackage{amsfonts} 
\usepackage{amssymb}
\usepackage{mathrsfs}
\usepackage{authblk}
\usepackage{cite}
\usepackage{typearea}                
\usepackage{todonotes}
\usepackage{scrlayer-scrpage} 
\usepackage{tikz,pgffor}
\usetikzlibrary{shadows}
\usepackage{mathtools}
\usepackage{tabularx}
\usepackage{rotating,graphics, graphicx, wrapfig}
\usepackage{pgffor}
\usepackage{a4wide}
\usepackage{graphicx}
\graphicspath{ {figures/} }
\usepackage{makeidx}
\usepackage[ruled]{algorithm2e} 
\usepackage{ gensymb }
\usepackage[makeroom]{cancel}
\makeindex

\usetikzlibrary{fit,shapes}
\usetikzlibrary{calc}
\usetikzlibrary{positioning}
\usetikzlibrary{shapes,arrows}
\usetikzlibrary{decorations.markings,decorations.pathmorphing}
\usepackage{xcolor}
\usepackage{fancybox}
\usepackage{xcolor,colortbl}
\usepackage{bibentry}
\usepackage{ marvosym }
\usepackage[xcolor,tikz]{mdframed}
\usepackage{rotating} 
\usepackage{mathtools}
\usetikzlibrary{decorations.pathreplacing}
\usepackage{array} 
\usepackage{rotating}
\usetikzlibrary{decorations.pathreplacing}
\usetikzlibrary{topaths,calc}
\usepackage{float}
\usepackage{cancel}
\usetikzlibrary{snakes}
\usetikzlibrary{tikzmark}

\ohead{\headmark}
\automark[subsection]{section}

\usetikzlibrary{decorations.pathreplacing}
\usetikzlibrary{topaths,calc}
\usepackage{float}
\usetikzlibrary{snakes,arrows,patterns}
\usetikzlibrary{tikzmark}

\tikzset{square arrow/.style={to path={-- ++(0,-.25) -| (\tikztotarget)}}}

\usetikzlibrary{positioning,calc}
\usetikzlibrary{arrows}
\usetikzlibrary{shapes.misc}

\newtheorem{defsatzusw}{}[section]
\newtheorem{definition}[defsatzusw]{Definition}
\newtheorem{theorem}[defsatzusw]{Theorem}

\newtheorem{lemma}[defsatzusw]{Lemma}
\newtheorem{corollary}[defsatzusw]{Corollary}

\newtheorem{remark}[defsatzusw]{Remark}

\newcommand{\qed}{\vspace{-.55cm}\begin{flushright}$\Box$\end{flushright}}

\newcommand{\N}{\mathbb N}

\newenvironment{proof}{\paragraph{$\mathit{Proof.}$}
	\vspace{-.4cm}\hspace{-.3cm}}{\vspace{-.55cm}\begin{flushright}$\Box$\end{flushright}}

\newenvironment{subproof}{\paragraph{}
	\vspace{-.8cm}\hspace{-.55cm}}{\vspace{-.75cm}\begin{flushright}$\blacksquare$\end{flushright}}

\definecolor{rwth-blue}{RGB}{0,84,159}
\definecolor{rwth-lblue}{RGB}{142,186,229}
\definecolor{rwth-lgreen}{RGB}{189,205,0}
\definecolor{rwth-green}{RGB}{87,171,39}
\definecolor{rwth-orange}{RGB}{246,168,0}
\definecolor{rwth-tuerkis}{RGB}{0,152,161}
\definecolor{rwth-red}{RGB}{204,7,30}
\definecolor{rwth-dred}{RGB}{161,16,53}

\definecolor{mygreen}{HTML}{086A08}
\definecolor{myred}{HTML}{9F0000}
\definecolor{myblue}{HTML}{08088A}
\definecolor{mygray}{HTML}{E6E6E6}

\pgfdeclarepatternformonly{newcrosshatch}{\pgfqpoint{-1pt}{-1pt}}{\pgfqpoint{4pt}{4pt}}{\pgfqpoint{3pt}{3pt}}%
{
	\pgfsetstrokeopacity{0.7}
	\pgfsetlinewidth{0.3pt}
	\pgfpathmoveto{\pgfqpoint{3.2pt}{0pt}}
	\pgfpathlineto{\pgfqpoint{0pt}{3.2pt}}
	\pgfpathmoveto{\pgfqpoint{0pt}{0pt}}
	\pgfpathlineto{\pgfqpoint{3.2pt}{3.2pt}}
	\pgfusepath{stroke}
}

\title{On the Domination Order among\\ Elimination Sequences}
\author[]{Michaela Hiller}
\affil[]{\small Lehrstuhl II f\"ur Mathematik, RWTH Aachen, 52062 Aachen, Germany }
\date{\vspace{-5ex}}

\begin{document}
	\setkomafont{sectioning}{\normalcolor\bfseries}	
	\maketitle

	\begin{abstract}\noindent
		\textbf{Abstract.}
		In 1991, it was shown by Favaron, Mahéo, and Saclé that the residue, which is defined as the number of zeros remaining when the Havel-Hakimi algorithm is applied to a degree sequence, 
		yields a lower bound on the independence number of any graph realising the sequence.
		In 1996, Triesch simplified and generalised the result by introducing elimination sequences.
		It was proved in \cite{KleitmanWang1973} that for any graphic sequence all elimination algorithms, i.e. laying-off vertices in any order, preserve that the sequence is graphic and terminate in a sequence of zeros.
		
		\noindent
		We now prove that for any degree sequence, the elimination sequence derived from the Havel-Havel algorithm dominates all other elimination sequences.
		Our result implies a conjecture posed by Michael Barrus in 2010: 
		When iteratively laying off degrees from a graphic sequence until only a list of zeros remains, the number of zeros is at most the residue of this sequence \cite{url}.
	\end{abstract}

	
	\section{Introduction}
	
%
%

Let $G=(V,E)$ be a pair of finite sets, where $E\subseteq\binom{V}{2}$. 
Then $G$ is called a finite, simple   {graph}\index{graph} with   {vertex set} $V$ and {edge set} $E$.
The {order}\index{order} of $G$ is defined as $|V|$ and denoted by $n$.
The {degree sequence}\index{degree!sequence} of a given graph $G$ is defined as the decreasingly ordered sequence of the degrees in $G$, denoted by $\pi(G)=\pi=(d_1\geq\ldots\geq d_n)$.
By $n_k$ we denote the number of degrees $k$, by $n_{\geq k}$ the number of degrees greater or equal to $k$ and by $n_{>k}$ the number of degrees greater than $k$.
We refer to a sequence $\pi$ as {graphic}\index{sequence!graphic} if there is a graph $G$ with $\pi(G)=\pi$. 
The {length} of a sequence $\pi$ is denoted by $\ell(\pi)$.
To merge two degree sequences $\pi_1$ and $\pi_2$, we write $(\pi_1,\pi_2)$ and assume that the degrees are resorted in descending order.
Analogously to an empty graph, we call a sequence of zeros {empty sequence}\index{sequence!empty}. 
Note that the empty sequence of any length $\ell(\pi)$ is graphic since it is realised by the empty graph of order $n=\ell(\pi)$.

Given a degree sequence $\pi=(d_1\geq\ldots\geq d_n)$, we define its {conjugate}\index{sequence!conjugate} as
$\overline{\pi}=(\overline{d_1},\dots \overline{d_{d_1}})$ where $\overline{d_i}=\big|\{j~|~d_j\geq i\}\big|$ for $i\in\{1,\dots,d_1\}$.
Partitions are frequently presented using Ferrers diagrams\index{Ferrers diagram}. 
The Ferrers diagram representing a degree sequence $\pi=(d_1\geq\ldots\geq d_n)$ consists of $n$ rows with $d_i$ dots in row $i$.
Conjugating a partition translates into flipping the Ferrers diagram diagonally.

To compare sequences, we introduce the domination\index{domination order} order: A sequence $\pi$ is said to {dominate} a sequence $\sigma$ if 
\begin{align*}
	&\sum\limits_{i=1}^{\ell(\pi)}\pi(i)=\sum\limits_{i=1}^{\ell(\sigma)}\sigma(i)\quad&&\text{and}\quad &&\sum\limits_{i=1}^p\pi(i)\geq\sum\limits_{i=1}^p\sigma(i)\quad&&\text{ for all }1\leq p\leq{\ell(\pi)},
\end{align*}
which we denote by $\pi\succeq\sigma$.

Following \cite{Triesch1996JoGT}, let $\pi=(d_1\geq\ldots\geq d_n)$, $d_{n+1}\coloneqq 0$, and $1\leq i<j \leq n+1$ such that $d_i\geq d_j+2$, where $i=\max\{p~|~d_p=d_i\}$ and $j=\min\{p~|~d_p=d_j\}$. 
When decreasing $d_i$ and increasing $d_j$ by one, we obtain a new sequence $\tau$ and say that $\tau$ was generated by an $(i,j)${-step}\index{$(i,j)$-step} from $\pi$. 
If $j=\min\{p~|~d_p\leq d_i-2\}$, the $(i,j)$-step is called elementary.
Then $\pi\succeq  \tau$.\\
In \cite{aigner1984}, it was proved that $\pi\succeq \rho$ if and only if $\rho$ can be obtained by a sequence of elementary $(i,j)$-steps from $\pi$.
From $\pi\succeq \rho$ and $\pi$ graphic, it follows that $\rho$ is also graphic \cite{RuchGutman1979}.\\
Inverting an $(i,j)$-step, i.e. increasing $d_i$ and decreasing $d_j$ by one, is called an {inverse} $(i,j)$-step. 
Clearly, the sequence obtained by an inverse $(i,j)$-step dominates the original sequence.\\

Havel \cite{Havel1955} and Hakimi \cite{hakimi1963realizability} characterised graphic sequences by  independently developing the same idea.

\begin{definition}[Havel-Hakimi Operator]
	Let $\pi=(d_1\geq\ldots\geq d_n)$ be a sequence where $n\geq d_1+1$. 
	The three-step of deleting term $d_1$, reducing each of the following $d_1$ terms by one, and resorting the obtained sequence in descending order is defined as the Havel-Hakimi operator -- also referred to as the Havel-Hakimi reduction step. We denote the application of this operator by
	$$\mathcal{H(\pi)}\coloneqq\left(d_2-1,\ldots,d_{d_1+1}-1,d_{d_1+2},\ldots,d_n\right).$$\index{Havel-Hakimi operator}
	\vspace{-.4cm}
\end{definition}
Note that in the notation above the terms are not yet reordered.

\begin{theorem}[Havel-Hakimi Theorem \textnormal{\cite{Havel1955,hakimi1963realizability}}]\label{Havel-Hakimi-Theorem}
	The sequence $\pi=(d_1\geq\ldots\geq d_n)$ is graphic if and only if $n\geq d_1+1$ and the reduced sequence $\mathcal{H(\pi)}$ is graphic.
\end{theorem}

\noindent
Repeated application of the Havel-Hakimi reduction step eventually leads to a sequence of $n-i$ zeros, i.e. $\mathcal{H}^{(i)}=(0^{n-i})$ for an $i\in\{0,1,\dots,n-1\}$.
We define the minimum $i$ such that $\mathcal{H}^i(\pi)=(0^{n-i})$ as {depth} $dp(\pi)$ of $\pi$ and the number of zeros as {residue}\index{residue} $R(\pi)=n-dp(\pi)$ of $\pi$.
For a graph $G$ with degree sequence $\pi$, the residue is given as $R(G)=R(\pi)$.\\
Note that if $\pi$ is a graphic sequence, then (with slight abuse of notation) $\pi':=(\pi,0^t)$  is graphic as well and $R(\pi')=R(\pi)+t$ for $t\in\N$.

Kleitman and Wang generalised the  definition and theorem stated above by allowing for an arbitrary vertex or degree to be laid off \cite{KleitmanWang1973}: \\
Given a sequence $\pi=(d_1\geq \ldots\geq d_n)$, we can choose any vertex $v_i$, remove the corresponding degree $d_i$ in $\pi$ and decrease the largest (other) $d_i$ degrees by one.
This corresponds to connecting the chosen vertex to the $d_i$ remaining other vertices with highest labels. 
Then $\pi$ is a graphic sequence if and only if the sequence we obtain by laying off $d_i$ is graphic.\\

The computer programme \textsc{Graffiti} proposed that the residue of a degree sequence bounds the independence number of all possible realisations \cite{Fajtlowicz1988}.
The claim was first proved by Favaron, Mahéo, and Saclé \cite{Favaron1991} in 1991, then simplified by Griggs and Kleitman \cite{GriggsKleitman1994} and generalised by Triesch \cite{Triesch1996JoGT}.

\begin{theorem}\label{residue_bound}
	Let $G$ be a simple graph. 
	Then $\alpha(G)\geq R(G)$.
\end{theorem}

\noindent
The proof given by Favaron, Mahéo, and Saclé \cite{Favaron1991} is premised on the residue preserving the dominance order for fixed $n$.

\begin{lemma}\label{residue_dominate}\textnormal{\cite{Favaron1991}}
	Given two degree sequences $\pi,\sigma$ with $\pi\succeq\sigma$ and $\ell(\pi)=\ell(\sigma)$, then $$R(\pi)\geq R(\sigma).$$
\end{lemma}

\noindent
For a shorter proof of the generalised version of Lemma \ref{residue_dominate}, Triesch defined the elimination sequence \textnormal{\cite{Triesch1996JoGT}}.

\begin{definition}[Elimination Sequence]
	Given a degree sequence $\pi=(d_1\geq\ldots\geq d_n)$ with depth $dp(\pi)=:dp$, the sequence 
	$${E}(\pi)\coloneqq\Big( \max(\pi), \max\left(\mathcal{H}(\pi)\right), \max\left(\mathcal{H}^2(\pi)\right), \dots, \max\left(\mathcal{H}^{dp-1}(\pi)\right), \underbrace{0, \dots,0}_{n-dp} \Big)$$
	is called the elimination sequence of $\pi$.\index{elimination sequence}
\end{definition}

\noindent
The main result in \cite{Triesch1996JoGT} is the preservation of the domination order when transitioning from degree sequences to the corresponding elimination sequences. 

\begin{theorem}\label{Dominanz_E}\textnormal{\cite{Triesch1996JoGT}}
	If $\pi$ and $\rho$ are two degree sequences with $\pi\succeq \sigma$, then $E(\pi)\succeq E(\sigma)$.
\end{theorem}

\noindent
We can now state our main theorem.

\begin{theorem}\label{domination_elimination_sequences}
	For a degree sequence $\pi$, the Havel-Hakimi elimination sequence ${E}(\pi)$ dominates any other elimination sequence $\tilde{E}(\pi)$.
\end{theorem}

\section{Preliminary}

To prove the above stated theorem, we need some preliminary lemmas.
\begin{lemma}\label{zusätzlicherGraderhältDominanz}
	Let $\pi$ and $\sigma$ be two sequences with $\pi\succeq \sigma$, then $(\pi, x)\succeq (\sigma, x)$ for all $x\in\N_0$.
\end{lemma}
\begin{proof}
	Let $\pi\coloneqq(\pi_1\geq\ldots\geq \pi_n)$ and $\sigma\coloneqq(\sigma_1\geq\ldots\geq \sigma_n)$ with $\pi\succeq \sigma$. Thus, for the conjugate sequences, it follows that $\overline{\pi}\preceq\overline{\sigma}$. 
	Now, $\overline{(\pi,x)}=(\overline{\pi}_1+1,\dots,\overline{\pi}_x+1,\overline{\pi}_{x+1},\dots,\overline{\pi}_{\pi_1})$ and $\overline{(\sigma,x)}=(\overline{\sigma}_1+1,\dots,\overline{\sigma}_x+1,\overline{\sigma}_{x+1},\dots,\overline{\sigma}_{\sigma_1})$.
	Therefore, $\overline{(\pi,x)}\preceq\overline{(\sigma,x)}$, which implies $(\pi, x)\succeq (\sigma, x)$.\vspace{-.1cm}
\end{proof}


\begin{lemma}\label{additiveDominanz}
	Suppose $\pi_1,\sigma_1,\pi_2,\sigma_2$ are graphic sequences with $\pi_1\succeq \sigma_1$ and $\pi_2\succeq \sigma_2$. Then $(\pi_1,\pi_2)\succeq(\sigma_1,\sigma_2)$.
\end{lemma}
\begin{proof}
	Repeated application of Lemma \ref{zusätzlicherGraderhältDominanz} yields $(\pi_1,\sigma_1)\succeq(\pi_1,\sigma_2)\succeq(\pi_2,\sigma_2)$ as claimed.\vspace{-.1cm}
\end{proof}

\noindent
Given a (not necessarily graphic) sequence $\pi$ of non-negative elements and $p,q\in\N_0$, $p\leq \ell(\pi)$, $q\leq \ell(\pi)$, let $\pi^p$ denote the sequence obtained after reducing the largest $p$ entries of $\pi$ by one and resorting the elements in descending order. 
Further, let $\pi^{p,q}$ denote the sequence we obtain by reducing the $q$ greatest entries in $\pi^p$ by one and again resorting the elements in descending order.

\begin{lemma}\label{reduce_a_a+1}
	Let $\pi$ be sequence with elements sorted in decreasing order and $p\in\N_0$, $p\leq n$ such that all entries in  $\pi^{p,p+1}$ and $\pi^{p+1,p}$ are still non-negative. Then $\pi^{p,p+1}=\pi^{p+1,p}$.
\end{lemma}

\begin{proof}
	Let the element at position $p+1$ in $ \pi$ be $k$, see Figure \ref{XXXX}.
	After the first reduction step, in $\pi^p$ and $\pi^{p+1}$ the two successive elements $k$ and $k-1$ occur with different frequencies, because in $\pi^{p+1}$ one (additional) entry $k$ was reduced by one to $k-1$, while in $\pi^{p}$ the element $k$ from position $p+1$ remains unchanged. 
	Therefore, in $\pi^{p+1}$, element $k$ occurs once less than in $\pi^{p}$, but element $k-1$ occurs once more than in $\pi^{p}$.
	Thus, $n_k(\pi^{p})=n_k(\pi^{p+1})+1$ and $n_{k-1}(\pi^{p})+1=n_{k-1}(\pi^{p+1})$.\\
	Further, the entry at position $p+1$ in both sequences, $\pi^{p}$ and $\pi^{p+1}$, is either $k$ or $k-1$: 
	The entry at position $p+1$ cannot be larger than $k$ since the element at position $p+1$ was $k$ in the first step and the entries have only be decreased after that.
	Also, the entry at position $p+1$ cannot be smaller than $k-1$ as the largest $p+1$ elements of $\pi$ are still greater than or equal to $k-1$ in $\pi^{p}$ and $\pi^{p+1}$. \\
	Now, assume that the element at position $p+1$ in $\pi^{p}$ is $k$.
	Then the largest $p$ entries in $\pi$ must have been greater or equal to $k+1$.
	Consequently, the element at position $p$ in $\pi^{p+1}$ is $k$ as well, due to rearrangement.
	When applying the second step to both sequences, we reduce one additional entry $k$ in $\pi^{p}$ compared to $\pi^{p+1}$. 
	Since the frequency of $k$ was one less and of $k-1$ one greater in $\pi^{p+1}$ than in $\pi^{p}$ after the first step, the frequencies are now balanced and $\pi^{p,p+1}=\pi^{p+1,p}$.\\
	If, however, the element at position $p+1$ in $\pi^{p}$ is $k-1$, the number of entry $k-1$ among the largest $p+1$ elements in $\pi^p$ equals the number of entry $k$ among the largest $p$ elements in $\pi$.
	This implies that we have the same number of entry $k-1$ among the first $p$ elements in $\pi^{p+1}$.
	Therefore, after the second step, $n_{k-2}(\pi^{p,p+1})=n_{k-2}(\pi^{p+1,p})$. 
	Furthermore, in $\pi^p$ we reduce one additional entry $k$ to $k-1$ compared to $\pi^{p+1}$, while the number of elements greater than $k$ are equal in $\pi^{p}$ and $\pi^{p+1}$. 
	Hence, $n_{k}(\pi^{p})=n_{k}(\pi^{p+1})+1$ and $n_{k-1}(\pi^{p})+1=n_{k-1}(\pi^{p+1})$ imply $n_{k}(\pi^{p,p+1})=n_{k}(\pi^{p+1,p})$ and $n_{k-1}(\pi^{p,p+1})=n_{k-1}(\pi^{p+1,p})$.
	Therefore, $\pi^{p,p+1}=\pi^{p+1,p}$, which was to be shown.\vspace{-.15cm}
\end{proof}

\begin{figure}[H]
	\begin{center}
		\begin{tikzpicture}[decoration=brace, node distance=5em, every node/.style={scale=.85}, scale=.85]
			\node(5) at (-1.4,-.7)[]{\textcolor{gray}{\tiny$p-2$}};
			\node(5) at (-.35,-.7)[]{\textcolor{gray}{\tiny$p-1$}};
			\node(5) at (0.3,-.7)[]{\textcolor{gray}{\tiny$p$}};
			\node(5) at (.94,-.7)[]{\textcolor{gray}{\tiny$p+1$}};
			\node(5) at (1.7,-.7)[]{\textcolor{gray}{\tiny$p+2$}};
			\node(5) at (2.7,-.7)[]{\textcolor{gray}{\tiny$p+2$}};
			\node(5) at (0,0)[]{\large$\pi=(\quad\ldots\quad,~k+1,~k,~k,~k,~k,~k-1,\quad\ldots\quad)$};
			\draw[-stealth] (-3,-1) -- (.62,-1) node[midway,sloped,below] {$p$};
			\draw[-stealth] (-3,-1.7) -- (1.25,-1.7) node[midway,sloped,below] {$p+1$};
			
		\end{tikzpicture}
		\caption{Reducing the largest $p$ and $p+1$ terms in $\pi$ to obtain $\pi^p$ and $\pi^{p+1}$, respectively.}
		\label{XXXX}		
	\end{center}
\end{figure}	

\begin{remark}\label{reduce_a_a+1_often}
	Let $\sigma$ be a (not necessarily graphic) sequence  of non-negative elements and $\sigma^{x_1,\dots,x_l}$ denote the sequence we obtain when successively reducing the largest $x_j$ elements in the sequence for $j\in\{1,\dots,l\}$ and resorting the entries in decreasing order after each step.
	Then iterative application of Lemma \ref{reduce_a_a+1} implies $\sigma^{p,(p+1)^x}=\sigma^{p+1^x,p}$.
\end{remark}

\noindent
For a stronger result, we prove the following generalisation.

\begin{lemma}\label{^a,b=^b,a}
	Assume $\pi=(d_1\geq\ldots\geq d_n)$ to be a degree sequence and  $p,q\in\N$, $p\leq n$, $q\leq n$ such that all entries in $\pi^{p,q}$ and $\pi^{q,p}$ are non-negative. Then $\pi^{p,q}=\pi^{q,p}$.
\end{lemma}

\begin{proof}
	Let $p<q$ with $d_p=k$ and $d_q=l$. Thus, $k\geq l$. 
	We define $\tilde{p}:=p-n_{>k}$ and $\tilde{q}:=q-n_{>k}$. 
	With this notation, it suffices to show that $\sigma^{\tilde{p},\tilde{q}}=\sigma^{\tilde{q},\tilde{p}}$ for $\sigma=(d_{n_{>k}+1},\dots,d_n)$ of length $n-n_{>k}$ since all degrees greater than $k$ are identically reduced by two in both $\pi^{p,q}$ and $\pi^{q,p}$.
	We consider four cases.
	\begin{enumerate}
		\item Case: $k=l$. \\
		If $\tilde{p}+\tilde{q}\leq n_k$, it follows that $n_k(\sigma^{\tilde{p},\tilde{q}})=n_k-\tilde{p}-\tilde{q}=n_k(\sigma^{\tilde{q},\tilde{p}})$ and $n_{k-1}(\sigma^{\tilde{p},\tilde{q}})=n_{k-1}+\tilde{p}+\tilde{q}=n_{k-1}(\sigma^{\tilde{q},\tilde{p}})$, while all other $n_i$ with $i\in\{k-2,\dots,1\}$ remain unchanged. 
		Therefore, $\sigma^{\tilde{p},\tilde{q}}=\sigma^{\tilde{q},\tilde{p}}$.\\
		If otherwise $\tilde{p}+\tilde{q}> n_k$, we obtain $n_k(\sigma^{\tilde{p}})=n_k-\tilde{p}$ and $n_k(\sigma^{\tilde{q}})=n_k-\tilde{q}$. 
		In both sequences, the remaining $n_k-\tilde{p}$ or $n_k-\tilde{q}$ entries of degree $k$ are also reduced after the second step. 
		Thus, we additionally decrease $\tilde{q}-(n_k-\tilde{p})$ entries of degree $k-1$ in $\sigma^{\tilde{p},\tilde{q}}$ and $\tilde{p}-(n_k-\tilde{q})$ in $\sigma^{\tilde{p},\tilde{q}}$. 
		Therefore, $\sigma^{\tilde{p},\tilde{q}}=\sigma^{\tilde{q},\tilde{p}}$.
		\begin{figure}[H]
			\begin{center}
				\begin{tikzpicture}
							\node(0) at (0,0-.2)[draw,  ellipse, minimum width=1pt]{\footnotesize$\sigma^{\tilde{p},\tilde{q}}$};
							\node(0) at (2.5,0-.2)[draw,  ellipse, minimum width=1pt]{\footnotesize$\sigma^{\tilde{q},\tilde{p}}$};
							\node(0) at (-1.5,-1)[]{\footnotesize \textcolor{gray}{$k$} };
							\node(0) at (-1.5,-1-1)[]{\footnotesize \textcolor{gray}{$k-1$} };
							\node(0) at (-1.5,-1-1-1)[]{\footnotesize \textcolor{gray}{$k-2$} };
							
							\draw (-.5,-1) -- (-.5,-1.25) -| (.5,-1);
							\node(1) at (0,-1)[]{};
							\draw (-.5,-1-1) -- (-.5,-1.25-1) -| (.5,-1-1);
							\node(11) at (0,-1-1.25)[]{};
							\draw (-.5,-1-1-1) -- (-.5,-1.25-1-1) -| (.5,-1-1-1);
							\node(111) at (0,-1-1.25-1)[]{};
							
							\draw (2.5-.5,-1) -- (2.5-.5,-1.25) -| (2.5+.5,-1);
							\node(2) at (2.5,-1)[]{};
							\draw (2.5-.5,-1-1) -- (2.5-.5,-1.25-1) -| (2.5+.5,-1-1);
							\node(22) at (2.5,-1-1.25)[]{};
							\draw (2.5-.5,-1-1-1) -- (2.5-.5,-1.25-1-1) -| (2.5+.5,-1-1-1);
							\node(222) at (2.5,-1-1.25-1)[]{};
							
							\draw[draw=orange, thick, -stealth] (1) to[bend right] node[left] {\textcolor{orange}{\scriptsize$\tilde{p}$}} (11);
							\draw[draw=rwth-blue, thick, -stealth] (1) to[bend left] node[right] {\textcolor{rwth-blue}{\scriptsize$n_k-\tilde{p}$}} (11);
							\draw[draw=rwth-blue, thick, -stealth] (11)+(.05,.05) to[bend left] node[right] {\textcolor{rwth-blue}{\scriptsize$\tilde{q}-(n_k-\tilde{p})$}} (111);
							
							\draw[draw=orange, thick, -stealth] (2) to[bend right] node[left] {\textcolor{orange}{\scriptsize$\tilde{q}$}} (22);
							\draw[draw=rwth-blue, thick, -stealth] (2) to[bend left] node[right] {\textcolor{rwth-blue}{\scriptsize$n_k-\tilde{q}$}} (22);
							\draw[draw=rwth-blue, thick, -stealth] (22)+(.05,.05) to[bend left] node[right] {\textcolor{rwth-blue}{\scriptsize$\tilde{p}-(n_k-\tilde{q})$}} (222);
				\end{tikzpicture}
			\caption{Reducing degrees in the case $k=l$ and $\tilde{p}+\tilde{q}> n_k$.}
			\end{center}
		\end{figure}
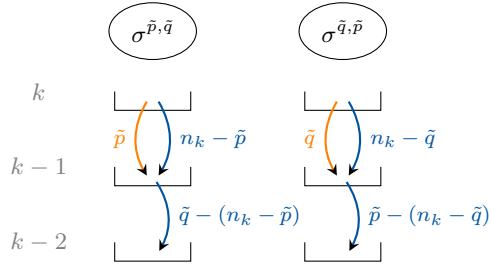
		\item Case: $k=l+1$.\\
		It suffices to consider $n_k$, $n_{k-1}$, and $n_{k-2}$ for $\sigma^{\tilde{p},\tilde{q}}$ and $\sigma^{\tilde{q},\tilde{p}}$ since all other degrees occur with equal frequencies in both sequences.\\
		In $\sigma^{\tilde{p}}$, we have $n_k(\sigma^{\tilde{p}})=n_k-\tilde{p}$, $n_{k-1}(\sigma^{\tilde{p}})=n_{k-1}+\tilde{p}$, and $n_{k-2}(\sigma^{\tilde{p}})=n_{k-2}$. 
		Hence, 
		\begin{align*}
			n_{k}(\sigma^{\tilde{p},\tilde{q}})&=n_k-\tilde{p}-(n_k-\tilde{p})=0,\\
			n_{k-1}(\sigma^{\tilde{p},\tilde{q}})&=n_{k-1}+\tilde{p}+(n_k-\tilde{p})-(\tilde{q}-(n_k-\tilde{p}))=n_{k-1}+2n_k-\tilde{p}-\tilde{q},\\
			n_{k-2}(\sigma^{\tilde{p},\tilde{q}})&=n_{k-2}+(\tilde{q}-(n_k-\tilde{p}))=n_{k-2}-n_k+\tilde{p}+\tilde{q}.
			\intertext{Now, in $\sigma^{\tilde{q}}$ we have $n_k(\sigma^{\tilde{q}})=0$, $n_{k-1}(\sigma^{\tilde{q}})=n_{k-1}+n_k-(\tilde{q}-n_k)$, and $n_{k-2}(\sigma^{\tilde{q}})=n_{k-2}+\tilde{q}-n_k$. 
				Thus, }
			n_{k}(\sigma^{\tilde{q},\tilde{p}})&=0,\\
			n_{k-1}(\sigma^{\tilde{q},\tilde{p}})&=n_{k-1}+n_k-(\tilde{q}-n_k)-\tilde{p}=n_{k-1}+2n_k-\tilde{p}-\tilde{q},\\
			n_{k-2}(\sigma^{\tilde{q},\tilde{p}})&=n_{k-2}+\tilde{q}-n_k+\tilde{p}=n_{k-2}-n_k+\tilde{p}+\tilde{q}.
		\end{align*}
		Again, $\sigma^{\tilde{p},\tilde{q}}=\sigma^{\tilde{q},\tilde{p}}$.
		\item Case: $k=l+2$.\\
		Analogous to the previous cases, we only need to compare $n_k$, $n_{k-1}$, $n_{k-2}$, and $n_{k-3}$ for $\sigma^{\tilde{p},\tilde{q}}$ and $\sigma^{\tilde{q},\tilde{p}}$.
		We obtain $n_k(\sigma^{\tilde{p}})=n_k-\tilde{p}$, $n_{k-1}(\sigma^{\tilde{p}})=n_{k-1}+\tilde{p}$, and $n_{k-2}(\sigma^{\tilde{p}})=n_{k-2}$, $n_{k-3}(\sigma^{\tilde{p}})=n_{k-3}$. Hence,  
		\begin{align*}
			n_{k}(\sigma^{\tilde{p},\tilde{q}})&=n_k-\tilde{p}-(n_k-\tilde{p})=0,\\
			n_{k-1}(\sigma^{\tilde{p},\tilde{q}})&=n_{k-1}+\tilde{p}+(n_k-\tilde{p})-(n_{k-1}+\tilde{p})=n_k-\tilde{p},\\
			n_{k-2}(\sigma^{\tilde{p},\tilde{q}})&=n_{k-2}+n_{k-1}+\tilde{p}-(\tilde{q}-n_k-n_{k-1})=n_{k-2}+2n_{k-1}+n_k+\tilde{p}-\tilde{q},\\
			n_{k-3}(\sigma^{\tilde{p},\tilde{q}})&=n_{k-3}+\tilde{q}-n_k-n_{k-1}.
			\intertext{Further, we have }
			n_{k}(\sigma^{\tilde{q}})&=n_k-n_k=0=n_{k}(\sigma^{\tilde{q},\tilde{p}}),\\
			n_{k-1}(\sigma^{\tilde{q}})&=n_{k-1}+n_k-n_{k-1}=n_k,\\
			n_{k-2}(\sigma^{\tilde{q}})&=n_{k-2}+n_{k-1}-(\tilde{q}-n_k-n_{k-1})=n_{k-2}+2n_{k-1}+n_k-\tilde{q},\\
			n_{k-3}(\sigma^{\tilde{q}})&=n_{k-3}+\tilde{q}-n_k-n_{k-1}=n_{k-3}(\sigma^{\tilde{q},\tilde{p}}),
		\end{align*}
		and therefore, $n_{k-1}(\sigma^{\tilde{q},\tilde{p}})=n_k-\tilde{p}$ and $n_{k-2}(\sigma^{\tilde{q},\tilde{p}})=n_{k-2}+2n_{k-1}+n_k-\tilde{q}+\tilde{p}$, which implies $\sigma^{\tilde{p},\tilde{q}}=\sigma^{\tilde{q},\tilde{p}}$.
		\item Case: $k>l+2$.\\
		In all other cases, we consider $n_k,n_{k-1},n_{k-2}$, and $n_l,n_{l-1}$ as, for $k-2>i>l$ and $l-1>i$, it clearly holds that $n_{i}(\sigma^{\tilde{p},\tilde{q}})=n_{i}(\sigma^{\tilde{q},\tilde{p}})$.
		Now, after the first step we obtain
		\begin{align*}
			n_k(\sigma^{\tilde{p}})&=n_k-\tilde{p},				&n_k(\sigma^{\tilde{q}})&=0,\\
			n_{k-1}(\sigma^{\tilde{p}})&=n_{k-1}+\tilde{p},		&n_{k-1}(\sigma^{\tilde{q}})&=n_{k},\\[4mm]
			n_{l}(\sigma^{\tilde{p}})&=n_l,						&n_{l}(\sigma^{\tilde{q}})&=n_{l}-\left(\tilde{q}-\sum_{i=l+1}^k n_i\right),\\
			n_{l-1}(\sigma^{\tilde{p}})&=n_{l-1},				&n_{l-1}(\sigma^{\tilde{q}})&=n_{l-1}+\tilde{q}-\sum_{i=l+1}^k n_i.
		\end{align*}
		Note that $n_{i}(\sigma^{\tilde{p}})=n_i$ and $n_{i}(\sigma^{\tilde{q}})=n_{i+1}$ for $k-1>i>l$. After the second step
		\begin{align*}
			n_k(\sigma^{\tilde{p},\tilde{q}})&=0,												&n_k(\sigma^{\tilde{q},\tilde{p}})&=0,\\
			n_{k-1}(\sigma^{\tilde{p},\tilde{q}})&=n_k(\sigma^{\tilde{p}})=n_k-\tilde{p},		&n_{k-1}(\sigma^{\tilde{q},\tilde{p}})&=n_k-\tilde{p},\\
			n_{k-2}(\sigma^{\tilde{p},\tilde{q}})&=n_{k-1}(\sigma^{\tilde{p}})=n_{k-1}+\tilde{p},\qquad	&n_{k-2}(\sigma^{\tilde{q},\tilde{p}})&=n_{k-1}+\tilde{p},\\[4mm]
			n_{l}(\sigma^{\tilde{p},\tilde{q}})&=n_{l}-\left(\tilde{q}-\sum_{i=l+1}^k n_i\right),										&n_{l}(\sigma^{\tilde{q},\tilde{p}})&=n_{l}(\sigma^{\tilde{q}})=n_{l}-\left(\tilde{q}-\sum_{i=l+1}^k n_i\right),\\
			n_{l-1}(\sigma^{\tilde{p},\tilde{q}})&=n_{l-1}+\tilde{q}-\sum_{i=l+1}^k n_i,										&n_{l-1}(\sigma^{\tilde{q},\tilde{p}})&=n_{l-1}(\sigma^{\tilde{q}})=n_{l-1}+\tilde{q}-\sum_{i=l+1}^k n_i,
		\end{align*}
		and therefore, $\sigma^{\tilde{p},\tilde{q}}=\sigma^{\tilde{q},\tilde{p}}$.\vspace{-.55cm}	
	\end{enumerate}
\end{proof}

\noindent
Given a degree sequence $\pi=(d_1\geq\ldots\geq d_n)$, let $\pi_i$ be the resulting sequence after laying off $d_i$.
Furthermore, let $\pi_{i,j}$ be the resulting sequence after first laying off $d_i$ and subsequently $d_j'$ where $\pi_i:=(d_1',\ldots, d'_{i-1},d'_{i+1}, d_n')$, according to Kleitman and Wang \cite{KleitmanWang1973}. 
If not all degrees of value $d_k$ are reduced when lying off $d_k$, we assume that $d_k$ itself remains unreduced.

\begin{lemma}\label{einfacherFall1}	
	Let $\pi=(d_1\geq\ldots\geq d_n)$ be a graphic sequence with $d_{1}\geq n_{\geq d_i}-1$ and $d_i\geq n_{d_1}$. Then $\pi_{1,i}=\pi_{i,1}$.
\end{lemma}

\begin{proof}
	Since $d_{1}+1\geq n_{\geq d_i}$, all entries of degree $d_i$ are reduced by one in $\pi_1$.
	Thus, in the second step, we lay off $d_i-1$ and reduce the largest $d_i-1$ entries in $\pi_1$, but then we can define $\sigma:=(d_2\geq\ldots\geq d_{i+1}\geq d_{i-1}\geq\ldots\geq d_n)$ and consider  $\sigma^{d_1-1,d_i-1}=\pi_{1,i}$, equivalently.\\ 
	Further, $d_i\geq n_{d_1}$ implies that, when first laying off $d_i$, all entries of degree $d_1$ are decreased by one.
	Hence,  in the second step we lay off $d_1-1$ and reduce the largest $d_1-1$ entries in $\pi_i$ to obtain $\pi_{i,1}$.
	Again, we can equivalently consider $\sigma^{d_i-1,d_1-1}=\pi_{i,1}$.\\
	Now, Lemma \ref{^a,b=^b,a} implies $\pi_{1,i}=\pi_{i,1}$, and thus, the assertion follows.\vspace{-.15cm}
\end{proof}

\begin{lemma}\label{einfacherFall2}
	Let $\pi=(d_1\geq\ldots\geq d_n)$ be a graphic sequence with $d_{1}< n_{\geq d_i}-1$ and $d_i< n_{d_1}$. Then $\pi_{1,i}=\pi_{i,1}$.
\end{lemma}

\begin{proof}
	As $d_{1}+1< n_{\geq d_i}$ implies that, when removing $d_1$, at least one entry of degree $d_i$ was not affected and is only deleted in the second step to obtain $\pi_{1,i}$, we can equivalently consider $\sigma:=(d_2\geq\ldots\geq d_{i+1}\geq d_{i-1}\geq\ldots\geq d_n)$ with $\sigma^{d_1,d_i}=\pi_{1,i}$.\\
	Also, from $d_i< n_{d_1}$ it follows that, when first laying off $d_i$, at least one entry of degree $d_1$ remains in $\pi_i$, which is then removed in the second step for $\pi_{i,1}$. 
	Again, we can equivalently consider $\sigma^{d_i,d_1}=\pi_{i,1}$.\\
	By Lemma \ref{^a,b=^b,a}, we obtain $\pi_{1,i}=\pi_{i,1}$, and the claim follows.\vspace{-.15cm}
\end{proof}


\section{Proof of the Main Theorem}

We can now show that the Havel-Hakimi elimination sequence dominates all other elimination sequences by an induction on the length of the sequence.
Our approach is based on the following idea:
By induction hypothesis the Havel-Hakimi elimination sequence $E$ dominates any elimination sequence $E'$ we obtain when first laying off a maximum degree (as in the first step of the Havel-Hakimi algorithm) and then following an arbitrary order. 
Further, let elimination sequence $\tilde{E}'$ result from first laying off an arbitrary degree $d_i$ and then following the Havel-Hakimi algorithm for the remaining sequence.
By $\tilde{E}$, we denote the elimination sequence we obtain when first laying off degree $d_i$ and then following any order. 
Then, also by the induction hypothesis, $\tilde{E}'$ dominates $\tilde{E}$.
Therefore, it remains to prove that $E'$ dominates $\tilde{E}'$.
$$E\stackrel{\text{(I.H.)}}{\succeq} E'\stackrel{}{\succeq} \tilde{E}'\stackrel{\text{(I.H.)}}{\succeq}\tilde{E}$$

%

\noindent
\textit{Proof of Theorem \ref{domination_elimination_sequences}.}
Let $R$ and $\tilde{R}$ denote the Havel-Hakimi residue and the residue of the modified algorithm, respectively.

For each elimination algorithm $\tilde{\mathcal{H}}\neq\mathcal{H}$, we can determine the first step in which $\tilde{\mathcal{H}}$ deviates from $\mathcal{H}$, while in all previous steps, the elimination elements equal those that occur in $E$.
Thus, by Lemma \ref{additiveDominanz} it suffices to show the claim for algorithms, where the first step differs from Havel-Hakimi, i.e. given a degree sequence $\pi=(d_1\geq\ldots\geq d_n)$, instead of $d_1$ we remove a degree $d_i<d_1$.

We proceed by induction on $n\in\N_0$, the length of the graphic sequence, and assume that an elimination algorithm does not stop until the sequence is empty, i.e. in the last $R$ or $\tilde{R}$ steps zeros are deleted but no entries are reduced. 
For the base case $n=0$, we clearly have $E(\pi)\succeq\tilde{E}(\pi)$.\\
By the induction hypothesis, $E(\pi)\succeq\tilde{E}(\pi)$ for all $\pi$ with at most $n$ degrees. 
Thus, we now assume $\pi$ to have $n+1$ entries.\\
First, consider the special case that $\pi$ requires less than two steps in $\mathcal{H}$ to be reduced to a sequence of zeros.
Then either $\pi$ is itself already a sequence of zeros and the claim immediately follows, or one Havel-Hakimi step is needed to obtain a null sequence.
Hence, we consider the last step in $\mathcal{H}$ before having a sequence of zeros, i.e. the  $(n-R-1)$-th step.
Therefore, 
$$\mathcal{H}^{(n-R-1)}(\pi)=(x,\underbrace{1,\dots,1}_{x\text{ times}},\underbrace{0,\dots,0}_{y\text{ times}})=:\sigma$$
for a certain $x\in\N$ with $y=R-x$, and $E(\sigma)=(x,0^R)$.
Now, removing an entry of degree one is the only possibility for an algorithm $\tilde{\mathcal{H}}\neq\mathcal{H}$ (apart from deleting a zero, which can be neglected because it has no effect on the subsequent steps or the elimination sequences).
From this it follows that 
$$\tilde{\mathcal{H}}(\sigma)=(x-1,\underbrace{1,\dots,1}_{x-1\text{ times}},\underbrace{0,\dots,0}_{y\text{ times}}).$$
In each subsequent step of $\tilde{\mathcal{H}}$, another entry of degree one is deleted, and in the final step (before obtaining a sequence of zeros) the current maximum degree is deleted. 
Consequently, in the penultimate step $j\geq1$, 
$$\tilde{\mathcal{H}}^{(j)}(\sigma)=(x-j,\underbrace{1,\dots,1}_{x-j\text{ times}},\underbrace{0,\dots,0}_{y\text{ times}}),$$
and $\tilde{E}(\sigma)=(x-j,\underbrace{1,\dots,1}_{j\text{ times}},\underbrace{0,\dots,0}_{R-j\text{ times}}).$
Hence, $E(\sigma)\succeq\tilde{E}(\sigma)$.\\

In the general case, where at least two Havel-Hakimi steps are needed to obtain a sequence of zeros, we suppose that $\tilde{\mathcal{H}}$ eliminates $d_i<d_1$ in the first step for any $d_i>0$. 
We can further consider, without loss of generality, the algorithm $\tilde{\mathcal{H}}'$ that after laying off $d_i$ always eliminates the current largest degree from the second step onwards since then, by the induction hypothesis, the elimination sequence of $\tilde{\mathcal{H}}'$ dominates the one of all other elimination algorithms $\tilde{\mathcal{H}}$, which first lay off $d_i$ as well but then follow any order.\\
In the following, let $d_j^\mathcal{\hat H}$ denote degree $d_j$ after the first elimination step of an elimination algorithm ${\mathcal{\hat H}}$. 
As mentioned above, if not all entries of degree $d_j$ are reduced in the laying-off step of ${\mathcal{\hat H}}$, we assume $d_j^{\mathcal{\hat H}}$ to remain degree $d_j$ if not otherwise stated. So, the following is proved for $d_i^{\mathcal{H}}=d_i$ in some cases and for $d_i^{\mathcal{H}}=d_i-1$ in others, both implying the claim. \\
For the sake of readability, we refrain noting the elements of the sequences in descending order and find
$$E(\pi)=(d_1,E(\pi_1))\succeq(d_1,d_i^\mathcal{H},E(\pi_{1,i}))\succeq E'(\pi)$$ 
and
$$\tilde{E}'(\pi)=(d_i,E(\pi_{i}))=(d_i,d_1^{\tilde{\mathcal{H}'}},E(\pi_{i,1}))\succeq \tilde{E}(\pi)$$ 
by the induction hypothesis and with slight abuse of notation. 
As seen above, it now remains to prove that $E'(\pi)\succeq\tilde{E}'(\pi)$.
In the following, we will show the stronger statement $(d_1,d_i^\mathcal{H},E(\pi_{1,i}))\succeq (d_i,d_1^{\tilde{\mathcal{H}'}},E(\pi_{i,1}))$.\\
In the case where $d_{1}\geq n_{\geq d_i}-1$ and $d_i\geq n_{d_1}$, Lemma \ref{einfacherFall1} implies $\pi_{1,i}=\pi_{i,1}$, and thus, $E(\pi_{1,i})=E(\pi_{i,1})$. 
Since $d_1^{\tilde{\mathcal{H}'}}+1=d_1>d_i=d_i^\mathcal{H}+1$, we obtain $(d_1,d_i^{\mathcal{H}}) \succeq (d_i,d_1^{\tilde{\mathcal{H}'}})$, and hence, $(d_1,d_i^\mathcal{H},E(\pi_{1,i}))\succeq (d_i,d_1^{\tilde{\mathcal{H}'}},E(\pi_{i,1}))$.\\
Analogously, in the case where $d_{1}< n_{\geq d_i}-1$ and $d_i< n_{d_1}$, it follows  by Lemma \ref{einfacherFall2}  that $\pi_{1,i}=\pi_{i,1}$, and thus, $E(\pi_{1,i})=E(\pi_{i,1})$.
As $d_1=d_1^{\tilde{\mathcal{H}'}}$ and $d_i=d_i^\mathcal{H}$, we obtain $(d_1,d_i^\mathcal{H},E(\pi_{1,i}))= (d_i,d_1^{\tilde{\mathcal{H}'}},E(\pi_{i,1}))$.\\
For the remaining cases, we have to prove that $E'\stackrel{}{\succeq} \tilde{E}'$ or, as mentioned above, the stronger statement $(d_1,d_i^\mathcal{H},E(\pi_{1,i}))\succeq (d_i,d_1^{\tilde{\mathcal{H}'}},E(\pi_{i,1}))$.
If not otherwise stated, we assume, without loss of generality, that $\mathcal{H'}$ lays off $d_1$, $d_i$, and corresponds to the Havel-Hakimi algorithm thereafter.

\begin{enumerate}
	\item Case: $d_{1}< n_{\geq d_i}-1$ and $d_i\geq n_{d_1}$.\\
	In this case, by laying off $d_1$, not all entries of degree $d_i$ are reduced in $\mathcal{H}'$, but, when laying off $d_i$ in the first step, all entries of degree $d_1$ are reduced in $\tilde{\mathcal{H}'}$.
	We distinguish two subcases.
	\begin{enumerate}
		\item[a)] $d_{1}> n_{> d_i}-1$:\\
		This means that removing and reducing $d_1$ in $\mathcal{H}'$ affects some degree $d_i$, which is then decreased by one.
		Unlike in the previous cases, we identify $d_i^{\mathcal{H}'}$ with this reduced degree.
		Analogously to the case where $d_{1}\geq n_{\geq d_i}-1$ and $d_i\geq n_{d_1}$ in Lemma \ref{einfacherFall1}, let $\sigma:=(d_2\geq\ldots\geq d_{i+1}\geq d_{i-1}\geq\ldots\geq d_n)$.
		Now,
		$\mathcal{H}'^{(2)}(\pi)=\pi_{1,i}=\sigma^{d_1-1,d_i-1}$.
		Furthermore, 
		$\tilde{\mathcal{H}}'^{(2)}(\pi)=\pi_{i,1}=\sigma^{d_i-1,d_1-1}$.
		Then Lemma \ref{^a,b=^b,a} implies $\pi_{1,i}=\pi_{i,1}$, and thus, the assertion follows.
		With $d_1^{\tilde{\mathcal{H}'}}+1=d_1>d_i=d_i^{\mathcal{H}'}+1$, it follows that  $(d_1,d_i^{\mathcal{H}'},E(\pi_{1,i}))\succeq (d_i,d_1^{\tilde{\mathcal{H}'}},E(\pi_{i,1}))$ by Lemma \ref{additiveDominanz}.\\
		\begin{figure}[H]
				\begin{tabular}{rcrc|c|ccc} 
					\hspace{2cm}$\mathcal{H}'$: & $d_1$ & ~$\ldots$~~~  $\ldots$ & $d_1$ & $\ldots$ & $d_i$ & $\ldots$   & $d_i$\\[2mm]
					& \phantom{$d_1-1$} &  $d_1-1$~  $\ldots$ & $d_1-1$ & $\ldots$ & $d_i-1$ &   $\ldots$ & $d_i$\\
				\end{tabular}
				\begin{tikzpicture}[overlay, remember picture]
					\begin{scope}
						\node[circle,draw,minimum width=14.5pt] at (-8.4,0.46){}; 
						\node[ellipse,draw,minimum width=33pt,minimum height=17pt] at (-2.45,-.2){}; 
						\draw[-stealth] (-7.7,.8) -- (-.9,.8){};
						\draw[] (-7.7,-.6) -- (-3,-.6){};
						\draw[-stealth] (-1.8,-.6) -- (-1.2,-.6){};
						\node at (1.8,0.2){$\rightarrow~~ \sigma^{d_1-1,d_i-1}$~~~~};
					\end{scope}
				\end{tikzpicture}\\[10mm]
		\begin{tabular}{lcrc|c|ccc} 
					\hspace{2cm}$\tilde{\mathcal{H}}'$:& $d_1$ &  \phantom{..}~$\ldots$~~~  $\ldots$ & $d_1$ & $\ldots$ & $d_i$ & $\ldots$   & $d_i$\\[2mm]
					&$d_1-1$  &  ~~$\ldots$~~~  $\ldots$ & $d_1-1$ & $\ldots$ &\phantom{$d_i-1$} &   $\ldots$ & $d_i$\\
				\end{tabular}
				\begin{tikzpicture}[overlay, remember picture]
					\begin{scope}
						\node[circle,draw,minimum width=14.5pt] at (-2.45,0.43){}; 
						\node[ellipse,draw,minimum width=35pt,minimum height=17pt] at (-8.4,-.2){}; 
						\draw[] (-8.9,.8) -- (-3,.8){};
						\draw[-stealth] (-1.8,-.6) -- (-.9,-.6){};
						\draw[] (-7.7,-.6) -- (-3,-.6){};
						\draw[-stealth] (-1.8,.8) -- (-1.2,.8){};
						\node at (1.8,0.2){$\rightarrow~~ \sigma^{d_i-1,d_1-1}$~~~~};
					\end{scope}
		
				\end{tikzpicture}
					\vspace{.3cm}
				\caption{If $d_{1}< n_{\geq d_i}-1$ and $d_i\geq n_{d_1}$, we obtain $\pi_{1,i}=\sigma^{d_1-1,d_i-1}$ and  $\pi_{i,1}=\sigma^{d_i-1,d_1-1}$ for the subcase $d_{1}> n_{> d_i}-1$.}
		\end{figure}
		\item[b)] $d_{1}\leq n_{> d_i}-1$:\\
		This implies that neither the removal of $d_1$ nor the removal of $d_i$ reduces any entry of degree $d_i$ in $\mathcal{H}'$.
		After the first step in $\tilde{\mathcal{H}'}$, where $d_i$ is laid off, we now increase ${d_1}^{\tilde{\mathcal{H}'}}=d_1-1$ by one to $d_1$ again and reduce the so far untouched degree $d_{d_i+1}$ by one,
		i.e. we apply an inverse $(1,d_i+1)$-step 
		(note that we slightly abuse our notation here since $d_{d_i+1}$ might have changed its position in the sequence). 
		Therefore, for the resulting sequence $\rho$, which is not necessarily graphic, it follows that $\rho\succeq\pi_i $.
		By this inverse $(1,{d_i+1})$-step, we imitate the second step of ${\mathcal{H}'}$, where $d_1$ is already deleted, and thus, the next greater degree $d_{d_i+1}$ is reduced additionally.
		Thereby, we ensure that the second step of $\tilde{\mathcal{H}'}$ (now applied on $\rho$) equals the first in ${\mathcal{H}'}$, in which $d_1$ was laid off.
		It follows that $\rho_1=\sigma^{d_1,d_i}=\pi_{1,i}$ with $\sigma$ as defined above and $\rho_1$ is graphic.
		Thus, $(d_1,d_i,E(\pi_{1,i}))=(d_1,d_i,E(\rho_1))\succeq (d_1,d_i,E(\pi_{i,1}))$ since $\rho$ arises from $\pi_i$ by an inverse $(1,d_i+1)$-step and with Theorem \ref{Dominanz_E}.
		
	\end{enumerate}
	\item Case: $d_{1}\geq n_{\geq d_i}-1$ and $d_i< n_{d_1}$.\\
	In this case, when laying off $d_1$, all entries of degree $d_i$ are reduced in $\mathcal{H}'$, but laying off $d_i$ in the first step does not reduce all entries of degree $d_1$ in $\tilde{\mathcal{H}}'$.
	It can easily be seen that $n_{d_1}\geq2$. 
	Also, $d_1>d_i+1$ or otherwise, $n_{d_1}=d_1$ and $n_{d_i}=1$: 
	Assume that $d_1=d_i+1$. 
	Then $d_i<n_{d_1}$ implies $d_1-1<n_{d_1}$ and $d_1\geq n_{\geq d_i}-1$ implies $d_1>n_{d_1}-1$. 
	Thus, $d_1=n_{d_1}$, and since $n_{\geq d_i}=n_{d_1}+n_{d_i}$ in this case, it follows that $n_{d_i}=1$.
	In the special case $d_1=n_{d_1}$, which implies $n_{d_i}=1$, as just seen, it follows that $d_i\geq2$, because otherwise, $\pi=(d_1^{d_1},1)$ is not graphic.\\
	In the general and in the special case, we continue the same way and define $x:=n_{d_1}-d_i>0$. 
	Note that (1) $n_{d_1}\geq x+1$ and (2) $d_1-x-2\geq0$ since $d_1>n_{d_1}$ or $d_i\geq2$. 
	Moreover, we have $d_{1}\geq n_{d_1}$ in any case.
	
	\textbf{Claim:}
	${\mathcal{H}'}^{(x+2)}(\pi)=\tilde{\mathcal{H}'}^{(x+2)}(\pi)$.
	
	\begin{subproof}
		Let $\pi_{\leftarrow}$ denote the $n_{d_1}$ entries of maximum degree and $\pi_{\rightarrow}$ the remaining $n-n_{d_1}$ entries. 
		We consider $\pi_{\leftarrow}$ and $\pi_{\rightarrow}$ separately. 
		These partial sequence are not required to be graphic. 
		Let $\mathcal{\hat H}(\pi_{\leftarrow})$ and $\mathcal{\hat H}(\pi_{\rightarrow})$ be the corresponding segments after applying elimination algorithm  $\mathcal{\hat H}$ to the whole sequence $\pi$.\\
		Note that for $d_1>n_{d_1}$ exchanges of degrees in $\pi_{\leftarrow}$ and $\pi_{\rightarrow}$ do occur neither in $\mathcal{H}'$ nor in $\tilde{\mathcal{H}'}$ in the first $x+2$ steps since (3) $d_1-x-2\geq d_i-1$ and $d_{1}\geq n_{\geq d_i}-1$.
		After the first step in ${\mathcal{H}'}$ and after the second in $\tilde{\mathcal{H}'}$, respectively, all former entries of degree $d_i$ in $\pi_{\rightarrow}$ are either removed or reduced, and all entries greater than $d_i$ in $\pi_{\rightarrow}$ are decreased simultaneously to the degrees in $\pi_{\leftarrow}$. 
		For the special case $d_1=n_{d_1}$ after step $x+2$ exchanges might occur, as we will see later. 
		However, after the last step we merely have to check that $n_{d_l}\left({\mathcal{H}'}^{(x+2)}(\pi)\right)=n_{d_l}\left({\tilde{\mathcal{H}'}}^{(x+2)}(\pi)\right)$ for all degrees $d_l$ and possible exchanges between $\pi_{\leftarrow}$ and $\pi_{\rightarrow}$ are irrelevant.\\\vspace{-.2cm}
		
		Applying $\mathcal{H}'$ to $\pi_{\leftarrow}$:
		
		By (1), applying $x+2$ steps of ${\mathcal{H}'}$ to $\pi$ corresponds to removing the largest degree and reducing all degrees in $\pi_{\leftarrow}$ by one in each step, except the second, where no degree is removed, but $d_i-1$ degrees (since $d_i$ was decreased in the first step in $\mathcal{H(\pi)}$) are reduced before $\pi_{\leftarrow}$ is being resorted in descending order.\\
		Thus, after the first step $n_{d_1}-1$ degrees $d_1-1$ are left and after the second, ${\mathcal{H}'}^{(2)}(\pi_{\leftarrow})$ consists of $n_{d_1}-d_i=x$ entries $d_1-1$ and $d_i-1$ entries $d_1-2$.\\
		For $2\leq s\leq x$, we have $n_{d_1}-d_i-s+2=x-s+2$ elements of maximum degree $d_1-s+1$ and $d_i-1$ degrees $d_1-s$ in ${\mathcal{H}'}^{(s)}(\pi_{\leftarrow})$.
		Then after step $x+1$, there is $n_{d_1}-d_i-x+1=1$ entry of maximum degree $d_1-x$ left, followed by $d_i-1$ entries of degree $d_1-x-1$.\\
		Finally, in ${\mathcal{H}'}^{(x+2)}(\pi_{\leftarrow})$, the last entry of degree $d_1-n_{d_1}+d_i=d_1-x$ is removed and again all other degrees are reduced by one. We obtain $d_i-1$ degrees of the new maximum degree $d_1-n_{d_1}+d_i-2=d_1-x-2$.\\\vspace{-.2cm}
		
		Applying $\tilde{\mathcal{H}'}$ to $\pi_{\leftarrow}$:
		
		We consider the effects of applying $x+2$ steps of $\tilde{\mathcal{H}'}$  to $\pi$. 
		In the first step, $d_i$ degrees are reduced by one and the partial sequence is resorted in descending order. 
		In all other steps, an entry of maximum degree is removed, while all other entries are reduced by one.
		Therefore, $\tilde{\mathcal{H}'}(\pi_{\leftarrow})$ consists of $n_{d_1}-d_i=x$ untouched degrees $d_1$ and $d_i$ degrees $d_1-1$. 
		In the second step, a degree $d_1$ is removed and all other degrees are reduced by one. 
		Thus, in $\tilde{\mathcal{H}'}^{(2)}(\pi_{\leftarrow})$ there are $x-1$ entries of degree $d_1-1$ and $d_i$ degrees $d_1-2$.\\
		Hence, for $2\leq s\leq x$, the sequence ${\mathcal{H}}'^{(s)}(\pi_{\leftarrow})$ has one additional entry of maximum degree, when compared to $\tilde{\mathcal{H}'}^{(s)}(\pi_{\leftarrow})$, and accordingly, one less of the second largest degree:
		In $\tilde{\mathcal{H}'}^{(s)}(\pi_{\leftarrow})$, there are $n_{d_1}-d_i-s+1=x-s+1$ degrees of maximum degree $d_1-s+1$ and $d_i$ degrees $d_1-s$.
		Then after step $x+1$, we obtain $n_{d_1}-d_i-(x+1)+1=0$ entries of degree $d_1-x$ and $d_i$ degrees of the new maximum degree $d_1-n_{d_1}+d_i-1=d_1-x-1$.\\
		Finally in $\tilde{\mathcal{H}'}^{(x+2)}(\pi_{\leftarrow})$, an entry of degree $d_1-n_{d_1}+d_i-1=d_1-x-1$ is deleted, and since $d_1-x-2\geq d_i-1$, all other entries are reduced by one, resulting in $d_i-1$ degrees $d_1-n_{d_1}+d_1-2=d_1-x-2$.
		This implies ${\mathcal{H}}'^{(x+2)}(\pi_{\leftarrow})=\tilde{\mathcal{H}'}^{(x+2)}(\pi_{\leftarrow})$.\\\vspace{-.2cm}
		
		Applying  $\mathcal{H}'$ and $\tilde{\mathcal{H}'}$ to $\pi_{\rightarrow}$:
		
		Next, we consider $\pi_{\rightarrow}$. We exclude $d_i$, because this degree is deleted in $\mathcal{H}'$ and $\tilde{\mathcal{H}'}$ in the second and first steps, respectively, and thus, merely affects $\pi_{\leftarrow}$. 
		However, we have to adapt the number of reduced degrees in the first step of $\mathcal{H}'(\pi_{\rightarrow})$, which are now $d_1-(n_{d_1}-1)-1=d_1-n_{d_1}=:p\in\N_0$. 
		In the second step of $\mathcal{H}'$ and in the first step of $\tilde{\mathcal{H}'}$, no degree in $\pi_{\rightarrow}$ is reduced.
		In the second step of $\tilde{\mathcal{H}'}$, we decrease $d_1-(n_{d_1}-1)=d_1-n_{d_1}+1=p+1$ degrees by one.
		For $2\leq s\leq x+1$, in the $s$-th step of both $\mathcal{H}'$ and $\tilde{\mathcal{H}'}$, we reduce $(d_1-s+1)-(n_{d_1}-s)=d_1-n_{d_1}+1=p+1$ entries of $\pi_{\rightarrow}$.
		Eventually, in the last step, i.e. the $(x+2)$-th step comparing the algorithms, there are $(d_1-x)-(d_i-1)=d_1-n_{d_1}+1=p+1$ degrees decreased in $\mathcal{H}'$ and $(d_1-x-1)-(d_i-1)=d_1-n_{d_1}=p$ in $\tilde{\mathcal{H}'}$. 
		Thus, it remains to prove that first reducing $x$ times $p+1$ entries and then once $p$, yields the same sequences as first reducing $p$ entries and then $x$ times $p+1$, where in both algorithms all degrees are resorted in descending order after each step.
		Remark \ref{reduce_a_a+1_often} yields exactly that, and thus, ${\mathcal{H}'}^{(x+2)}(\pi_{\rightarrow})=\tilde{\mathcal{H}'}^{(x+2)}(\pi_{\rightarrow})$.
		
		Therefore, we have ${\mathcal{H}'}^{(x+2)}(\pi)=\tilde{\mathcal{H}'}^{(x+2)}(\pi)$ as claimed.\vspace{-.25cm}
	\end{subproof}

	\textbf{Claim:}
	Let $E^\prime_{(x+2)}$ and $\tilde{E'}_{(x+2)}$ be the descendingly ordered sequences of the first $x+2$ eliminated degrees in $\mathcal{H}'$ and $\tilde{\mathcal{H}'}$, respectively. 
	Then $E^\prime_{(x+2)}\succeq\tilde{E'}_{(x+2)}$.
	
	\begin{subproof}
		The unsorted sequences of the first $x+2$ eliminated degrees of $\mathcal{H}'$ and $\tilde{\mathcal{H}'}$ are given by
		\begin{align*}
			\hspace{2.5cm}&(d_1,&&d_i-1,&&d_1-1,&& d_1-2,&&\dots,&&d_1-x\phantom{-1})\qquad\text{and}\hspace{15cm}\\
			&(d_i,&&d_1,&&d_1-1,&& d_1-2,&&\dots,&&d_1-x-1),
		\end{align*}
		respectively. Sorting the degrees in descending order yields
		\begin{align*}
			\hspace{.8cm}{E^\prime_{(x+2)}}=~&(d_1,&&d_1-1,&& d_1-2,&&\dots,&&d_1-x+1,&&d_1-x\phantom{-1},&& d_i-1)\qquad\text{and} \hspace{15cm}\\
			\tilde{E'}_{(x+2)}~=~&(d_1,&&d_1-1,&& d_1-2,&&\dots,&&d_1-x+1,&&d_1-x-1,&&d_i~~\phantom{-1}),
		\end{align*}  
		since $d_1-n_{d_1}\geq 1$ implies $d_1-x-1=d_1-n_{d_1}+d_i-1\geq d_i$. Thus, $E_{(x+2)}^\prime\succeq\tilde{E'}_{(x+2)}$.\\
	\end{subproof}
	Finally, Lemma \ref{additiveDominanz} implies $E'(\pi)=\left(E^\prime_{(x+2)},E(\tau)\right)\succeq\left(\tilde{E'}_{(x+2)},E(\tau)\right)=\tilde{E'}(\pi)$ for $\tau:={\mathcal{H}'}^{(x+2)}(\pi)=\tilde{\mathcal{H}'}^{(x+2)}(\pi)$. 
	This completes the proof of Theorem \ref{domination_elimination_sequences}.\vspace{-.425cm}\qed
\end{enumerate}

\noindent
Now, the conjecture by Barrus follows \cite{url}. 

\begin{corollary}
	For a degree sequence $\pi$, the residue $R(\pi)$ is always greater than or equal to the number of zeros $\tilde{R}(\pi)$ we obtain by laying off terms in any order deviating from the one imposed by the Havel-Hakimi algorithm.
\end{corollary}

\noindent
Note that the elimination sequences from different elimination algorithms can be incomparable.
Consider for example the graphic sequence $(4,3,3,3,2,1)$ and two elimination algorithms 
\begin{align*}
	\mathcal{H}_1:~&~~\bcancel{4}&&3&&3&&3&&2&&1\qquad\qquad &\mathcal{H}_2:~&~~4&&\bcancel{3}&&3&&3&&2&&1\\
	&&&2&&2&&2&&1&&\bcancel{1}		&&~~\bcancel{3}&&&&2&&2&&2&&1\\
	&&&2&&2&&1&&\bcancel{1}&&		&&&&&&\bcancel{1}&&1&&1&&1\\
	&&&2&&1&&\bcancel{1}&&&&		&&&&&&&&\bcancel{1}&&1&&0\\
	&&&1&&\bcancel{1}&&&&&&			&&&&&&&&&&0&&0\\
	&&&0&&&&&&&&					&&&&&&&&&&&&\\[2mm]	
	\rightarrow E_1=~&(4,&&1,&&1,&&1,&&1,&&0)	&\rightarrow E_2=~&(3,&&3,&&1,&&1,&&0,&&0)								 
\end{align*}
where for the corresponding elimination sequences we have $E_1\not\succeq E_2$ and $E_2\not\succeq E_1$.
Thus, the set of eliminations sequences is partially but not totally ordered.


	\bibliographystyle{abbrv} 
	\bibliography{Literatur_Gradsequenzen}

\begin{thebibliography}{10}

\bibitem{aigner1984}
M.~Aigner.
\newblock Uses of the diagram lattice.
\newblock {\em Mitt. Math. Semin. Gießen}, 163:61--77, 1984.

\bibitem{url}
M.~D. Barrus.
\newblock Havel-hakimi residue and independent sets.
\newblock {\em REGS 2010}, faculty.math.
  illinois.edu/$\sim$west/regs/hhresidue.html.

\bibitem{Fajtlowicz1988}
S.~Fajtlowicz.
\newblock On conjectures of graffiti.
\newblock {\em Annals of Discrete Mathematics}, 38:113--118, 1988.

\bibitem{Favaron1991}
O.~Favaron, M.~Mah{\'e}o, and J.-F. Sacl{\'e}.
\newblock On the residue of a graph.
\newblock {\em Journal of Graph Theory}, 15(1):39--64, 1991.

\bibitem{GriggsKleitman1994}
J.~R. Griggs and D.~J. Kleitman.
\newblock Independence and the {H}avel-{H}akimi residue.
\newblock {\em Discrete Mathematics}, 127(1-3):209--212, 1994.

\bibitem{hakimi1963realizability}
S.~L. Hakimi.
\newblock On realizability of a set of integers as degrees of the vertices of a
  linear graph {II}. uniqueness.
\newblock {\em Journal of the Society for Industrial and Applied Mathematics},
  11(1):135--147, 1963.

\bibitem{Havel1955}
V.~Havel.
\newblock A remark on the existence of finite graphs (hungarian).
\newblock {\em Casopis Pest. Mat.}, 80:477--480, 1955.

\bibitem{KleitmanWang1973}
D.~J. Kleitman and D.-L. Wang.
\newblock Algorithms for constructing graphs and digraphs with given valences
  and factors.
\newblock {\em Discrete Mathematics}, 6(1):79--88, 1973.

\bibitem{RuchGutman1979}
E.~Ruch and I.~Gutman.
\newblock The branching extent of graphs.
\newblock {\em J. Combin. Inform. System Sci}, 4(4):285--295, 1979.

\bibitem{Triesch1996JoGT}
E.~Triesch.
\newblock {Degree Sequences of Graphs and Dominance Order}.
\newblock {\em Journal of Graph Theory}, 22(1):89--93, 1996.

\end{thebibliography}
\end{document}